\documentclass[12pt]{article}
\usepackage{amsmath}
\usepackage{amssymb}
\renewcommand\smallskip{\vskip\smallskipamount}
\renewcommand\medskip{\vskip\medskipamount}
\renewcommand\bigskip{\vskip\bigskipamount}

\newcommand{\qed}{\hfill $\Box$ \medskip}

\begin{document}

\footnotetext{The author is partially supported by NSF Grant
DMS-0707086 and a Sloan Research Fellowship.}

\begin{center}
\begin{large}
\textbf{On the Local Solvability of Darboux's Equation}
\end{large}

\bigskip\medskip

MARCUS A. KHURI

\bigskip\medskip

\end{center}

\begin{abstract}
  We reduce the question of local nonsolvability of the Darboux
equation, and hence of the isometric embedding problem for
surfaces, to the local nonsolvability of a simple linear equation
whose type is explicitly determined by the Gaussian curvature.
\end{abstract}
\bigskip
\setcounter{equation}{0}

   Let $(M^{2},g)$ be a two-dimensional Riemannian manifold.  A
well-known problem is to ask, when can one realize this locally as
a small piece of a surface in $\mathbb{R}^{3}$?  That is, if the
metric $g=g_{ij}dx^{i}dx^{j}$ is given in the neighborhood of a
point, say $(x^{1},x^{2})=0$, when do there exist functions
$z_{\alpha}(x^{1},x^{2})$, $\alpha=1,2,3$, defined in a possibly
smaller domain such that $g=dz_{1}^{2}+dz_{2}^{2}+dz_{3}^{2}$?
This equation may be written in local coordinates as the following
determined system
\begin{equation*}
\sum_{\alpha=1}^{3}\frac{\partial z_{\alpha}}{\partial
x^{i}}\frac{\partial z_{\alpha}}{\partial x^{j}}=g_{ij}.
\end{equation*}
Due to its severe degeneracy, in the sense that every direction
happens to be a characteristic direction, little information has
been obtained by studying this system directly.  However a more
successful approach has been to reduce this system to the
following single equation of Monge-Amp\`{e}re type, known as the
Darboux equation:
\begin{equation}
\det\nabla_{ij}z=K|g|(1-|\nabla_{g}z|^{2})
\end{equation}
where $\nabla_{ij}$ are second covariant derivatives, $K$ is the
Gaussian curvature, $\nabla_{g}$ is the gradient with respect to
$g$, and $|g|=\det g$.  In fact, the local isometric embedding
problem is equivalent to the local solvability of this equation
(see the appendix).\par
   Let us first recall the known results.  Since equation (1) is
elliptic if $K>0$, hyperbolic if $K<0$, and of mixed type if $K$
changes sign, the manner in which $K$ vanishes will play the
primary role in the hypotheses  of any result.  The classical
results state that a solution always exists in the case that $g$
is analytic or $K(0)\neq 0$; these results may be found in [4].
C.-S. Lin provides an affirmative answer in [10] and [11] when $g$
is sufficiently smooth and satisfies $K\geq 0$, or $K(0)=0$ and
$\nabla K(0)\neq 0$.  When $K\leq 0$ and $\nabla K$ possesses a
certain nondegeneracy, Han, Hong, and Lin [5] show that a smooth
solution always exists if $g$ is smooth.  Lastly if the Gaussian
curvature vanishes to finite order and the zero set $K^{-1}(0)$
consists of Lipschitz curves intersecting tranvsversely, then Han
and the author [6] have proven the existence of smooth solutions
if $g$ is smooth.  Related results may be found in [1], [2], [3],
[7], [8].\par
   A negative result has been obtained by Pogorelov [13] (see also [12]),
who found a $C^{2,1}$ metric with no local $C^{2}$ isometric
embedding in $\mathbb{R}^{3}$.  More recently, the author [9] has
constructed $C^{\infty}$ examples of degenerate hyperbolic and
mixed type Monge-Amp\`{e}re equations of the form
\begin{equation}
\det (\partial_{ij}z+a_{ij}(p,z,\nabla z))=k(p,z,\nabla z)
\end{equation}
which do not admit a local solution, where $p=(x^{1},x^{2})$ and
$\partial_{ij}$ denote second partial derivatives.  A fundamental
part of the strategy in [9] is to reduce the local nonsolvability
of (2), to the local nonsolvability of a quasilinear equation
whose type is explicitly determined by the function $k$.  It is
the purpose of this article to show that the Darboux equation
possesses a similar property for a large class of Gaussian
curvatures.\par
   We begin by partially constructing the Gaussian curvature.
Here we will denote the coordinates $x^{1}$ and $x^{2}$ by $x$ and
$y$ respectively.  Define sequences of disjoint open squares
$\{X^{n}\}_{n=1}^{\infty}$ and $\{X_{1}^{n}\}_{n=1}^{\infty}$
whose sides are aligned with the $x$ and $y$-axes, and such that
$X^{n}$, and $X_{1}^{n}$ are centered at $q_{n}=(\frac{1}{n},0)$,
$X^{n}\subset X_{1}^{n}$, and $X^{n}$, $X_{1}^{n}$ have widths
$\frac{1}{2n(n+1)}$, $\frac{1}{n(n+1)}$, respectively.  Set
$K\equiv 0$ in $\mathbb{R}^{2}-\bigcup_{n=1}^{\infty}X_{1}^{n}$.
Define
\begin{equation*}
X=\{(x,y)\mid |x|<1, |y|<1\}
\end{equation*}
and let $\phi\in C^{\infty}(\overline{X})$ be such that $\phi$
vanishes to infinite order on $\partial X$, and either $\phi(q)>0$
or $\phi(q)<0$ for all $q\in X$ (here $\overline{X}$ denotes the
closure of $X$).  We now define $K$ in $X^{n}$ by
\begin{equation*}
K(q)=\gamma_{n}\phi(4n(n+1)(q-q_{n})),\text{ }\text{ }\text{
}\text{ }q\in \overline{X}^{n},
\end{equation*}
where $\{\gamma\}_{n=1}^{\infty}$ is a sequence of positive
numbers that are to be chosen with the property that
$\lim_{n\rightarrow\infty}\gamma_{n}=0$ in order to insure that
$K\in C^{\infty}(\mathbb{R}^{2})$.  A description of how $K$
should be prescribed in the remaining region
$\bigcup_{n=1}^{\infty}(X_{1}^{n}-X^{n})$ shall be given
below.\medskip

\textbf{Theorem 1.}  \textit{Suppose that $K$ adheres to the
description given above, and that a local $C^{5}$ solution $z$ of
the Darboux equation exists in a domain containing the origin.
Then in a neighborhood of a point on $\partial X^{n}$ for some $n$
sufficiently large, there exists a $C^{2}$ function $u$
constructed from $z$ which after an appropriate change of
coordinates satisfies the equation}
\begin{equation}
\partial_{tt}u+K\partial_{ss}u=Kf,
\end{equation}
\textit{where $f\in C^{0}$ also depends on $z$ and is strictly
positive.}\medskip

  This theorem suggests a strategy for constructing
smooth counterexamples to the local solvability of the Darboux
equation, or equivalently the local isometric embedding problem.
Namely, complete the construction of a smooth Gaussian curvature
function in the region $\bigcup_{n=1}^{\infty}(X_{1}^{n}-X^{n})$,
in such a way that the linear equation (3) can have no local
solution.  Whether this is possible is still an open question,
however as pointed out above, a similar strategy was successfully
employed for the related Monge-Amp\`{e}re equation (2). Note that
in order for this strategy to be utilized for the Darboux
equation, it must be shown that given a smooth function $K$ there
always exists a locally defined smooth metric $g$ having Gaussian
curvature $K$. This may be accomplished in the following way.  Let
$\Omega$ be a neighborhood of the origin, and let $G\in
C^{\infty}(\Omega)$ be the unique solution of the equation
\begin{equation*}
\partial_{xx}G+KG=0,\text{ }\text{ }\text{ }G(0,y)=1,\text{ }\text{
}\text{ }\partial_{x}G(0,y)=0.
\end{equation*}
By choosing $\Omega$ sufficiently small we have that $G>0$.  Then
\begin{equation*}
g=dx^{2}+G^{2}dy^{2}
\end{equation*}
is a smooth Riemannian metric and has Gaussian curvature $K$ in
the domain $\Omega$.\par
  The first step in verifying Theorem 1, will be to show that
certain second covariant derivatives of any solution of (1) cannot
vanish on $\partial X^{n}$ for $n$ sufficiently large. Suppose
that a local solution $z\in C^{2}$ of (1) exists, so that upon
rewriting the equation we have
\begin{equation}
b^{ij}\nabla_{ij}z=2K(1-|\nabla_{g}z|^{2}),
\end{equation}
where the Einstein summation convention concerning raised and
lowered indices has been used (this convention will also be
utilized in what follows) and
\begin{equation*}
(b^{ij})=|g|^{-1}\left(\begin{array}{cc}
\nabla_{22}z & -\nabla_{12}z \\
-\nabla_{12}z & \nabla_{11}z
\end{array}\right).
\end{equation*}
Then integrating by parts yields
\begin{equation}
\int_{X^{n}}b^{ij}\nabla_{ij}z
d\omega_{g}=-\int_{X^{n}}\nabla_{j}z\nabla_{i}b^{ij}d\omega_{g}
+\int_{\partial X^{n}}b^{ij}n_{i}\nabla_{j}z d\sigma_{g},
\end{equation}
where $d\omega_{g}$ and $d\sigma_{g}$ are the elements of area and
length with respect to $g$, and $(n_{1},n_{2})$ is the unit outer
normal to $\partial X^{n}$ also with respect to $g$.  In order to
calculate the interior term on the right-hand side we note that
$b^{ij}$ is a contravariant 2-tensor, so that
\begin{equation*}
\nabla_{i}b^{ij}=\partial_{i}b^{ij}+\Gamma_{il}^{i}b^{lj}+\Gamma_{il}^{j}b^{il}
\end{equation*}
where $\Gamma_{ij}^{l}$ are Christoffel symbols.  Therefore
\begin{eqnarray*}
\nabla_{i}b^{i1}&=&|g|^{-1}(\partial_{1}\nabla_{22}z-\partial_{2}\nabla_{12}z)
+|g|^{-2}(-\partial_{1}|g|\nabla_{22}z+\partial_{2}|g|\nabla_{12}z)\\
&
&+|g|^{-3/2}(\partial_{1}|g|^{1/2}\nabla_{22}z-\partial_{2}|g|^{1/2}\nabla_{12}z)
+\Gamma_{il}^{1}b^{il}\\
&=&|g|^{-1}(\partial_{1}\nabla_{22}z-\partial_{2}\nabla_{12}z+|g|\Gamma_{ij}^{1}b^{ij})
-\Gamma_{ij}^{i}b^{j1},
\end{eqnarray*}
after making use of the identity
\begin{equation*}
\Gamma_{ij}^{i}=|g|^{-1/2}\partial_{j}|g|^{1/2}.
\end{equation*}
Moreover direct computation shows that
\begin{eqnarray*}
&
&\partial_{1}\nabla_{22}z-\partial_{2}\nabla_{12}z+|g|\Gamma_{ij}^{1}b^{ij}\\
&=&
-\Gamma_{j2}^{j}\partial_{12}z+\Gamma_{j1}^{j}\partial_{22}z\\
&
&+(\partial_{2}\Gamma_{12}^{i}-\partial_{1}\Gamma_{22}^{i}-\Gamma_{11}^{1}\Gamma_{22}^{i}
+2\Gamma_{12}^{1}\Gamma_{12}^{i}-\Gamma_{22}^{1}\Gamma_{11}^{i})\partial_{i}z\\
&=&|g|(\Gamma_{j2}^{j}b^{12}+\Gamma_{j1}^{j}b^{11})\\
&
&+(\partial_{2}\Gamma_{12}^{i}-\partial_{1}\Gamma_{22}^{i}-\Gamma_{11}^{1}\Gamma_{22}^{i}
+2\Gamma_{12}^{1}\Gamma_{12}^{i}-\Gamma_{22}^{1}\Gamma_{11}^{i}-\Gamma_{j2}^{j}\Gamma_{12}^{i}
+\Gamma_{j1}^{j}\Gamma_{22}^{i})\partial_{i}z,
\end{eqnarray*}
and we observe that the coefficient of $\partial_{i}z$ is in fact
a curvature term.  More precisely, if it is denoted by $\chi^{i}$
then
\begin{equation*}
\chi^{i}=\partial_{2}\Gamma_{12}^{i}-\partial_{1}\Gamma_{22}^{i}+
\Gamma_{12}^{j}\Gamma_{j2}^{i}-\Gamma_{22}^{j}\Gamma_{j1}^{i}=-R_{212}^{i}
=-g^{i1}|g|K
\end{equation*}
where $R_{jkl}^{i}$ is the Riemann tensor.  We now have
\begin{equation*}
\partial_{1}\nabla_{22}z-\partial_{2}\nabla_{12}z+|g|\Gamma_{ij}^{1}b^{ij}
=|g|(\Gamma_{j2}^{j}b^{12}+\Gamma_{j1}^{j}b^{22}-g^{i1}K\partial_{i}z)
\end{equation*}
so that
\begin{equation}
\nabla_{i}b^{i1}=-Kz^{1}.
\end{equation}
Similarly
\begin{equation}
\nabla_{i}b^{i2}=-Kz^{2}.
\end{equation}
With the help of (4), (6), and (7) it follows that (5) becomes
\begin{eqnarray}
& &\int_{X^{n}}K(2-3|\nabla_{g}z|^{2})d\omega_{g}\\
&=&\int_{\partial X^{n}}
|g|^{-1/2}[(\nabla_{1}z\nabla_{22}z-\nabla_{2}z\nabla_{12}z)\overline{n}_{1}
+(\nabla_{2}z\nabla_{11}z-\nabla_{1}z\nabla_{12}z)\overline{n}_{2}]d\sigma,\nonumber
\end{eqnarray}
where $(\overline{n}_{1},\overline{n}_{2})$ is the Euclidean unit
outer normal to $\partial X^{n}$ and $d\sigma$ is the Euclidean
element of length.\par
  The integral equality (8) will now be used to show that certain
second covariant derivatives of any solution of the Darboux
equation cannot vanish on $\partial X^{n}$ for $n$ sufficiently
large. Let $-\mathrm{v}_{n}$, $+\mathrm{v}_{n}$ represent the left
and right vertical portions of $\partial X^{n}$, respectively, and
let $+\mathrm{h}_{n}$, $-\mathrm{h}_{n}$ represent the top and
bottom horizontal portions of $\partial X^{n}$,
respectively.\medskip

\textbf{Lemma 2.}  \textit{Suppose that $K$ satisfies the
hypotheses of Theorem 1.  Then it is not possible for a $C^{2}$
solution $z$ of (1) to satisfy the following property for any $n$
sufficiently large:}
\begin{equation}
\nabla_{22}z|_{\pm\mathrm{v}_{n}}=0,\text{ }\text{ }\text{ }\text{
}\nabla_{11}z|_{\pm\mathrm{h}_{n}}=0.
\end{equation}

\textit{Proof.}  We proceed by contradiction and assume that
property (9) holds.  Then since $K|_{\partial X^{n}}=0$, the
Darboux equation implies that $\nabla_{12}z|_{\partial X^{n}}=0$.
Therefore the right-hand side of (8) vanishes.  However this
yields a contradiction, as the left-hand side is nonzero for large
$n$.  To see this last fact observe that according to the
appendix, any solution of the Darboux equation yields an isometric
embedding $F=(z_{1},z_{2},z)$ of the metric $g$.  So that by
performing an appropriate rigid body motion of this embedding, to
obtain $\overline{F}=AF$ where $A$ is an orthogonal matrix, we can
ensure that the new third component $\overline{z}$ of
$\overline{F}$ satisfies $|\nabla\overline{z}|(0,0)=0$.
Furthermore the appendix also shows that $\overline{z}$ must
satisfy the Darboux equation, and so we have
$2-3|\nabla_{g}\overline{z}|^{2}>1$ inside $X^{n}$ if $n$ is
chosen sufficiently large.  Therefore since $K$ never vanishes on
$X^{n}$, integral equality (8) yields a contradiction. \qed

  In light of Lemma 2, there must exist a point $p\in\partial X^{n}$
at which one of the given second covariant derivatives is nonzero.
As arguments similar to those presented below may be applied if
$p\in-\mathrm{v}_{n}$ or $p\in\pm\mathrm{h}_{n}$, we assume
without loss of generality that $p\in+\mathrm{v}_{n}$ so that
$\nabla_{22}z(p)\neq 0$.  It follows that after a change of
coordinates near $p$, a solution $u$ of equation (3) may be
constructed.  The following lemma will complete the proof of
Theorem 1.\medskip

\textbf{Lemma 3.}  \textit{Suppose that there exists a $C^{5}$
solution $z$ of the Darboux equation satisfying
$\nabla_{22}z(p)\neq 0$.  Then there exists a $C^{3}$ local change
of coordinates near $p=(p^{1},p^{2})$ given by}
\begin{equation*}
t=x-p^{1},\text{ }\text{ }\text{ }\text{ }s=s(x,y),
\end{equation*}
\textit{and a $C^{2}$ solution $u$ of the equation}
\begin{equation*}
\partial_{tt}u+K\partial_{ss}u=Kf,
\end{equation*}
\textit{where $f\in C^{0}$ and is strictly positive if $n$ is
sufficiently large.}\medskip

\textit{Proof.}  The desired coordinates $(t,s)$ will be chosen to
eliminate the mixed second covariant derivative appearing in (4).
Since $b^{ij}$ is a contravariant 2-tensor, under a coordinate
change $\overline{x}^{i}=\overline{x}^{i}(x^{1},x^{2})$ it
transforms by
\begin{equation*}
\overline{b}^{ij}=b^{lm}\frac{\partial\overline{x}^{i}}{\partial
x^{l}}\frac{\partial\overline{x}^{j}}{\partial x^{m}}.
\end{equation*}
Therefore by setting $t=\overline{x}^{1}=x-p^{1}$, we seek
$s=\overline{x}^{2}$ such that
\begin{equation}
\overline{b}^{12}=b^{11}\partial_{x}s+b^{12}\partial_{y}s=0,\text{
}\text{ }\text{ }\text{ }s(p^{1},y)=cy,
\end{equation}
where $c$ is a nonzero constant to be determined.  Observe that
since $b^{11}=|g|^{-1}\nabla_{22}z\neq 0$ near $p$, the line
$x=p^{1}$ is noncharacteristic for (10).  Thus the theory of first
order partial differential equations guarantees the existence of a
unique local solution $s\in C^{3}$, in light of the fact that
$b^{11},b^{12}\in C^{3}$.\par
  We now calculate each of the new coefficients for the Darboux
equation.  First note that $\overline{b}^{11}=b^{11}$, and with
the help of (10)
\begin{eqnarray*}
\overline{b}^{22}&=&b^{11}(\partial_{x}s)^{2}+2b^{12}\partial_{x}s\partial_{y}s
+b^{22}(\partial_{y}s)^{2}\\
&=&(b^{11})^{-1}(\partial_{y}s)^{2}\det
b^{ij}\\
&=&(|g|b^{11})^{-1}(\partial_{y}s)^{2}K(1-|\nabla_{g}z|^{2}).
\end{eqnarray*}
Therefore in the new coordinates Darboux's equation (4) is given
by
\begin{equation}
b^{11}\overline{\nabla}_{11}z+K\overline{f}\overline{\nabla}_{22}z=2K(1-|\nabla_{g}z|^{2}),
\end{equation}
where $\overline{\nabla}_{ij}$ denote covariant derivatives with
respect to the new coordinates $(t,s)$ and
\begin{equation*}
\overline{f}=(|g|b^{11})^{-1}(\partial_{y}s)^{2}(1-|\nabla_{g}z|^{2}).
\end{equation*}
Notice that if we choose
\begin{equation*}
c=b^{11}|g|^{1/2}(1-|\nabla_{g}z|^{2})^{-1/2}(p),
\end{equation*}
then $(b^{11})^{-1}\overline{f}(p)=1$.  Moreover by setting
\begin{equation*}
u(t,s)=z(t,s)-\int_{0}^{t}\left(\int_{0}^{t^{'}}(\overline{\Gamma}_{11}^{1}\partial_{t}z
+\overline{\Gamma}_{11}^{2}\partial_{s}z)(t^{''},s)dt^{''}\right)dt^{'}
\end{equation*}
we have $\partial_{tt}u=\overline{\nabla}_{11}z$, so that (11)
becomes
\begin{equation*}
\partial_{tt}u+K\partial_{ss}u=Kf
\end{equation*}
with
\begin{equation*}
f=(b^{11})^{-1}[2(1-|\nabla_{g}z|^{2})+(\overline{f}(p)-\overline{f})\overline{\nabla}_{22}z]
+\overline{\Gamma}_{22}^{1}\partial_{t}z+\overline{\Gamma}_{22}^{2}\partial_{s}z
+\partial_{ss}(u-z).
\end{equation*}
Lastly we observe that $f(t,s)>0$ in a sufficiently small
neighborhood of $p$ if $n$ is large, since as in the proof of
Lemma 2 we may assume that $|\nabla z|(0,0)=0$.  \qed

\begin{center}
\textbf{Appendix}
\end{center}

   Here we show that the local isometric embedding problem is equivalent
to the local solvability of the Darboux equation (1).  Assume that
there exists a local $C^{2}$ embedding $F=(z_{1},z_{2},z_{3})$ for
a given metric $g$.  Then according to the Gauss equations
\begin{equation*}
\nabla_{ij}F=h_{ij}\nu,
\end{equation*}
where $h_{ij}$ are the components of the second fundamental form
with respect to a unit normal $\nu$.  Then by taking the Euclidean
inner product of this equation with the vector $\vec{k}=(0,0,1)$,
we obtain
\begin{equation*}
\det\nabla_{ij}z=K|g|(\nu\cdot \vec{k})^{2}
\end{equation*}
where for convenience we denote $z_{3}$ by $z$.  Furthermore, if
$\times$ represents the cross product operation between two
vectors in $\mathbb{R}^{3}$ then
\begin{equation*}
(\nu\cdot
\vec{k})^{2}=1-\left|\frac{(\partial_{1}F\times\partial_{2}F)\times
\vec{k}}{|\partial_{1}F\times
\partial_{2}F|}\right|^{2}=1-g^{ij}\partial_{i}z\partial_{j}z=1-|\nabla_{g}z|^{2},
\end{equation*}
where $g^{ij}$ are components of the inverse matrix
$(g_{ij})^{-1}$.  Clearly the remaining two components of $F$ must
also satisfy equation (1).  Conversely, if a local solution of (1)
exists for a given metric $g$ and $|\nabla_{g}z|<1$, then a
calculation shows that $g-dz^{2}$ is a Riemannian metric and is
flat. It follows that there exists a local change of coordinates
$z_{1}=z_{1}(x^{1},x^{2})$, $z_{2}=z_{1}(x^{1},x^{2})$ such that
$g-dz^{2}=dz_{1}^{2}+dz_{2}^{2}$.

\begin{center}\bigskip
\textbf{References}
\end{center}\bigskip


\noindent[1]\hspace{.06in} Q. Han, \textit{On the isometric
embedding of surfaces with Gauss curvature
changing}\par\hspace{.06in}\textit{sign cleanly}, Comm. Pure Appl.
Math., $\mathbf{58}$ (2005), 285-295, MR 2094852.\bigskip

\noindent[2]\hspace{.06in} Q. Han, \textit{Local isometric
embedding of surfaces with Gauss curvature
changing}\par\hspace{.06in}\textit{sign stably across a curve},
Calc. Var. \& P.D.E., $\mathbf{25}$ (2006), 79-103, MR
2183856.\bigskip

\noindent[3]\hspace{.06in} Q. Han, \textit{Smooth local isometric
embedding of surfaces with Gauss
curvature}\par\hspace{.06in}\textit{changing sign cleanly},
preprint.\bigskip

\noindent[4]\hspace{.06in} Q. Han, \& J.-X. Hong,
\textit{Isometric Embedding of Riemannian Manifolds in
Eu-}\par\hspace{.06in}\textit{clidean Spaces}, Mathematical
Surveys and Monographs, Vol. 130, AMS,
Provi-\par\hspace{.06in}dence, RI, 2006, MR 2261749.\bigskip

\noindent[5]\hspace{.06in} Q. Han, J.-X. Hong, \& C.-S. Lin,
\textit{Local isometric embedding of surfaces
with}\par\hspace{.06in}\textit{nonpositive Gaussian curvature,} J.
Differential Geom., $\mathbf{63}$ (2003), 475-520,
MR\par\hspace{.06in}2015470.\bigskip

\noindent[6]\hspace{.06in} Q. Han, \& M. Khuri, \textit{On the
local isometric embedding in $\mathbb{R}^{3}$ of surfaces
with}\par\hspace{.06in}\textit{Gaussian curvature of mixed sign},
preprint.\bigskip

\noindent[7]\hspace{.06in} M. Khuri, \textit{The local isometric
embedding in $\mathbb{R}^{3}$ of two-dimensional
Riemannian}\par\hspace{.06in}\textit{manifolds with Gaussian
curvature changing sign to finite order on a
curve},\par\hspace{.06in}J. Differential Geom., $\mathbf{76}$
(2007), 249-291, MR 2330415.\bigskip

\noindent[8]\hspace{.06in} M. Khuri, \textit{Local solvability of
degenerate Monge-Amp\`{e}re equations and
applica-}\par\hspace{.06in}\textit{tions to geometry}, Electron.
J. Diff. Eqns., $\mathbf{2007}$ (2007), No. 65, 1-37,
MR\par\hspace{.06in}2308865.\bigskip

\noindent[9]\hspace{.06in} M. Khuri, \textit{Counterexamples to
the local solvability of Monge-Amp\`{e}re
equations}\par\hspace{.06in}\textit{in the plane,} Comm. PDE,
$\mathbf{32}$ (2007), 665-674, MR 2334827.\bigskip

\noindent[10]  C.-S. Lin, \textit{The local isometric embedding in
$\mathbb{R}^{3}$ of} 2-\textit{dimensional
Riemannian}\par\hspace{.06in}\textit{manifolds with nonnegative
curvature}, J. Differential Geom., $\mathbf{21}$ (1985), no.
2,\par\hspace{.02in} 213-230, MR 0816670.\bigskip

\noindent[11]  C.-S. Lin, \textit{The local isometric embedding in
$\mathbb{R}^{3}$ of two-dimensional
Riemannian}\par\hspace{.06in}\textit{manifolds with Gaussian
curvature changing sign cleanly}, Comm. Pure
Appl.\par\hspace{.06in}Math., $\mathbf{39}$ (1986), no. 6,
867-887, MR 0859276.\bigskip






\noindent[12]  N. Nadirashvili, \& Y. Yuan, \textit{Improving
Pogorelov's isometric embedding coun-}\par \text{
}\textit{terexample}, Calc. Var. Partial Differential Equations,
$\mathbf{32}$ (2008), no. 3, 319-323,\par \text{ }MR
2393070.\bigskip

\noindent[13]  A. Pogorelov, \textit{An example of a
two-dimensional Riemannian metric not}\par\textit{ admitting a
local realization in $E_{3}$}, Dokl. Akad. Nauk. USSR,
$\mathbf{198}$ (1971),\par \text{ }42-43, MR 0286034.\smallskip



\bigskip\bigskip\footnotesize

\noindent\textsc{Department of Mathematics, Stony Brook
University, Stony Brook, NY 11794}\par

\noindent\textit{E-mail address}: \verb"khuri@math.sunysb.edu"

\end{document}